\documentclass[11pt,reqno]{amsart}
\usepackage[a4paper,margin=27mm]{geometry}
\usepackage[T1]{fontenc}
\usepackage{graphicx}
\usepackage[section]{placeins}
\usepackage{needspace}
\usepackage{lmodern,amsmath,amssymb,amsthm,mathtools,microtype,booktabs,mathrsfs}
\usepackage[hidelinks]{hyperref}
\hypersetup{pdftitle={Cubic local observability of Gaussian-blind sources in thermal phi-four fields},pdfauthor={Rubin Zhang and Ning Zhang}}
\numberwithin{equation}{section}
\newtheorem{mainthm}{Theorem}

\newtheorem{theorem}{Theorem}[section]
\newtheorem{lemma}[theorem]{Lemma}
\newtheorem{proposition}[theorem]{Proposition}
\newtheorem{corollary}[theorem]{Corollary}
\theoremstyle{definition}

\newtheorem{remark}[theorem]{Remark}
\newtheorem{problem}[theorem]{Problem}
\newcommand{\E}{\mathbb E}
\newcommand{\R}{\mathbb R}
\newcommand{\Z}{\mathbb Z}
\newcommand{\T}{\mathcal T}
\newcommand{\K}{\mathcal K_\beta}

\newcommand{\OO}{\mathcal F_O}
\newcommand{\dd}{\,\mathrm d}
\newcommand{\ip}[2]{\langle #1,#2\rangle}
\newcommand{\norm}[1]{\lVert #1\rVert}
\newcommand{\KL}{\operatorname{KL}}
\newcommand{\Var}{\operatorname{Var}}

\newcommand{\supp}{\operatorname{supp}}
\newcommand{\op}{\mathrm{op}}
\newcommand{\LL}{\mathfrak L}
\newcommand{\cA}{\mathfrak A}
\newcommand{\cF}{\mathfrak F}

\DeclareMathOperator{\Tr}{Tr}
\allowdisplaybreaks[2]
\title[Cubic observability of thermal sources]{Cubic local observability of Gaussian-blind sources in thermal $\phi^4_2$ fields}
\author{Rubin Zhang}
\address{School of Mathematics and Statistics, Huazhong University of Science and Technology, Wuhan, Hubei 430074, China}
\curraddr{Department of Mathematical and Statistical Sciences, University of Alberta, Edmonton, Alberta T6G 2G1, Canada}
\email{rbzhang3@163.com}
\author{Ning Zhang}
\address{School of Mathematics and Statistics, Huazhong University of Science and Technology, Wuhan, Hubei 430074, China}
\email{nzhang2@hust.edu.cn}
\thanks{Corresponding author: Ning Zhang (nzhang2@hust.edu.cn).}
\date{November 2026}
\begin{document}
\begin{abstract}
We study external sources that are invisible to a spatially local observation in a massive Gaussian thermal field on a circle. A quadratic-reference calibration removes the first-order thermal mass-shift signal. For the resulting quartic model we prove a cutoff-independent, source-uniform $C^1$ expansion of local Wick cubic means. Its leading operator has the time-integrated cube of the thermal covariance as its kernel. A positive bilateral Laplace representation proves injectivity on the entire external source space, while analyticity gives exponential singular-value upper bounds. The expansion yields local Araki relative-entropy lower bounds on every fixed finite-dimensional source space, and a uniform separation theorem at every prescribed Sobolev resolution. Under an additional interval-length condition, a reflection-positivity identity yields an explicit, nonoptimal double-logarithmic conditional stability bound for this specific thermal kernel. Continuum probability estimates and the integrability assumptions of the Feynman--Kac--Nelson construction are verified in Section~\ref{sec:construct}. Uniqueness at a single fixed positive coupling is formulated as an open problem (Problem~\ref{prob:fixed}).
\end{abstract}
\maketitle

\section{Introduction}
A static source changes both a thermal state and the information available to a local observer. In a free massive field the source Bogoliubov--Kubo--Mori (BKM) metric is the inverse massive Laplacian. This global identity does not imply local identifiability: an exterior source can translate the field only outside the observation region. We study the local information carried by a calibrated quartic response and establish separate linearized and finite-coupling stability statements.

Fix a spatial circle $\mathbb T_L$, an inverse temperature $\beta>0$, and a mass $m>0$. Let $I$ and $O$ be nonempty open intervals with disjoint closures. After choosing an origin and a lift to $\R$, write
\begin{equation}\label{eq:geom}
 I=(0,\ell),\qquad \overline O\subset(\ell-L,0),\qquad 0<\ell<L.
\end{equation}
For $A=-\partial_x^2+m^2$, every $v\in C_c^\infty(I)$ generates a Gaussian-blind source $Av$: the corresponding free mean is $v$, and its restriction to $O$ vanishes. In fact, the entire equal-time local Gaussian state is unchanged, not only its mean.

A necessary qualification concerns thermal mass calibration. Vacuum and thermal Wick ordering differ by a quadratic term. The resulting first-order mean signal can already be produced by a shifted Gaussian mass. We therefore retain an explicit quadratic coefficient $b\ge0$ and use the calibrated source
\begin{equation}\label{eq:Jcalintro}
 J_{\lambda,v}=(A+2b\lambda)v.
\end{equation}
The Gaussian comparison model with quadratic operator $A+2b\lambda$ is then blind to $J_{\lambda,v}$ for every $\lambda\ge0$. Thermal Wick ordering corresponds to $b=0$; the original vacuum-ordered quartic interaction corresponds, at this fixed temperature and up to a scalar, to $b=6d_\beta$, where $d_\beta$ is defined in \eqref{eq:dbeta}. This is an exact calibration of the specified quadratic reference, not an identification of the full interacting pole mass.

Let $\mu$ be the centered Gaussian measure on $\mathbb T_\beta\times\mathbb T_L$ with covariance $C=(-\partial_\tau^2+A)^{-1}$. Thermal Wick powers are denoted by colons. The local observable
\[
 W_h=\int_O h(y):\phi(0,y)^3:\dd y,\qquad h\in L^2(O),
\]
is a spatially smeared Wick cubic field, not the cube of a smeared linear field. Its mean in the calibrated source state defines a map $\mathcal M_\lambda:H_0^1(I)\to L^2(O)$ by duality. All quantum-state assertions below are made for smooth sources; the Euclidean observation map extends to all of $H_0^1(I)$.

\subsection{Principal conclusions}
The main statements separate the continuum cubic response, the properties of its linearized kernel, and the resulting local information bounds.

\begin{mainthm}[Uniform continuum cubic response; see Theorem~\ref{thm:response}]
Let $H=H_0^1(I)$ with $\norm v_H^2=\int_I(|v'|^2+m^2v^2)\dd x$ and $Y=L^2(O)$. For each $R<\infty$ and fixed $\lambda_*<\infty$, the calibrated observation map satisfies
\begin{equation}\label{eq:mainexpintro}
 \mathcal M_\lambda(v)=-24\lambda\T v+\mathcal R_\lambda(v),\qquad
 \sup_{\norm v_H\le R}\norm{D\mathcal R_\lambda(v)}_{H\to Y}\le C_R\lambda^2,
 \quad 0\le\lambda\le\lambda_*,
\end{equation}
where
\begin{equation}\label{eq:Tintro}
 (\T v)(y)=\int_I\K(x-y)v(x)\dd x,\qquad
 \K(z)=\int_0^\beta C(\tau,z)^3\dd\tau.
\end{equation}
The coefficient $24=4\cdot3!$ is independent of the quadratic calibration parameter $b$.

\end{mainthm}

\begin{mainthm}[Injectivity and analytic smoothing; see Theorem~\ref{thm:injective} and Proposition~\ref{prop:compact}]
The operator $\T:L^2(I)\to L^2(O)$ is injective, compact, and of infinite rank. Its singular values satisfy $s_n(\T)\le Ce^{-cn}$ for some $C,c>0$. The injectivity proof uses positive Laplace atoms and a half-plane zero-set argument, while analyticity yields the singular-value bound.

\end{mainthm}

\begin{mainthm}[Local information separation; see Theorems~\ref{thm:finite} and~\ref{thm:resolution}]
For every fixed finite-dimensional smooth source space $E\subset H$ and fixed source radius $R$, sufficiently small positive coupling gives
\[
 S_O(\omega_{\lambda,v}\Vert\omega_{\lambda,w})
 \ge c_{E,R}\lambda^2\norm{v-w}_H^2
 \quad(v,w\in E,\ \norm v_H,\norm w_H\le R),
\]
where $c_{E,R}>0$ depends on $E$, $R$, and the fixed model parameters. For a compact infinite-dimensional Dirichlet spectral prior and each prescribed resolution $r>0$, one sufficiently small positive coupling separates all source pairs at $H$-distance at least $r$. The admissible coupling $\lambda_r$ depends on $r$; for smooth pairs the separation also holds for local Araki relative entropy.

\end{mainthm}

\begin{mainthm}[Conditional double-logarithmic stability; see Theorem~\ref{thm:doublelog}]
If $\ell<L/2$, then for $s>0$ and $u\in X_s$ with $\norm u_{X_s}\le M$,
\begin{equation}\label{eq:rateintro}
 \norm u_H\le
 \frac{C_sM}{\LL(M/\norm{\T u}_{Y})^{s/(1+s)}},
 \qquad \LL(x)=1+\log(1+\log(1+x)).
\end{equation}
The rate follows from a reflected positive energy and a quantitative bound for the associated noninteger moment system. The space $X_s$ is the Dirichlet spectral prior defined in Section~\ref{sec:entropy}.

\end{mainthm}

\subsection{Relation to previous work}
The BKM metric and canonical correlations are classical \cite{Kubo,Petz}. Finite-temperature free-energy convexity also enters Hamiltonian learning \cite{AnshuEtAl}. Our setting concerns the information available in a local algebra from an exterior continuum source. The constructive $P(\phi)_2$ framework and the relation between periodic Euclidean processes and KMS states are developed in \cite{HoeghKrohn,Simon,GlimmJaffe,KleinLandau,GJcut,GJfull}. We use the positive-temperature perturbation theorem \cite[Theorem~7.11 and Propositions~7.8--7.10]{GJcut} after verifying its sharp-time hypotheses.

Local relative entropy for free-field coherent excitations was studied in \cite{Longo,Ciolli,Casini}. For interacting theories, Drago, Faldino and Pinamonti \cite{DragoFaldinoPinamonti} analyze perturbative equilibrium-state relative entropy, while Fr\"ob, Much and Papadopoulos \cite{FrobMuchPapadopoulos} compute a first-order quartic-coupling contribution for vacuum and coherent states in a Rindler wedge. In contrast, we obtain lower bounds for constructive finite-temperature states from a measured local Wick cubic observable and data processing \cite{Araki}, without computing an interacting modular Hamiltonian.

Nonlinear inverse-source methods have precedents in hyperbolic and elliptic problems \cite{KurylevLassasUhlmann,FeizmohammadiOksanen,LiimatainenLin}. The inverse scattering analysis of Sasaki and Suzuki \cite{SasakiSuzuki} concerns a linearly sourced Klein--Gordon field. The model-specific results here are injectivity of the calibrated thermal cubic kernel and local information separation from its continuum response.

Propagation of smallness for analytic functions and heat observability are studied in \cite{Vessella,Apraiz}; we use the spectral heat kernel on the circle \cite{Davies}. Analytic-kernel approximation and classical M\"untz theory supply ingredients for the singular-value and moment estimates \cite{Trefethen,BorweinErdelyi}. Related quantitative moment and inverse-Laplace methods appear in \cite{Saitoh,GarciaFerreroRuland,Gendron,Choulli}, and inverse-problem instability is illustrated in \cite{Mandache,Natterer}. Section~\ref{sec:stability} identifies the kernel-specific reflected energy and the source of the double-logarithmic loss. The final regularization step uses standard Hilbert-space arguments \cite{Engl}.

\section{Gaussian blindness and the calibrated interacting model}\label{sec:model}
\subsection{The two-moment obstruction}
The massive periodic Green function is
\begin{equation}\label{eq:green}
 G_\kappa(z)=\frac{e^{-\kappa z}+e^{-\kappa(L-z)}}{2\kappa(1-e^{-\kappa L})},\qquad 0<z<L.
\end{equation}
For a smooth source supported in $I$, its free mean on $O$ is
\begin{equation}\label{eq:twomom}
 (A^{-1}J)(y)=\frac{e^{my}\int_Ie^{-mx}J(x)\dd x
 +e^{-m(L+y)}\int_Ie^{mx}J(x)\dd x}{2m(1-e^{-mL})}.
\end{equation}
\begin{proposition}[Exact Gaussian blindness]\label{prop:blind}
On the smooth source space $C_c^\infty(I)$, the restriction of the free thermal state to the common equal-time local algebra on $O$ depends exactly on the two moments in \eqref{eq:twomom}. The kernel of this restriction, relative to the zero source, is
\[
 \{Av:v\in C_c^\infty(I)\}.
\]
In particular this kernel has codimension two in $C_c^\infty(I)$.
\end{proposition}
\begin{proof}
Completing the square gives a Weyl displacement of the centered free thermal state by the mean $A^{-1}J$. The field covariance is source independent and the momentum mean is zero. Hence local Weyl expectations depend only on the restriction of $A^{-1}J$ to $O$. The two exponentials in \eqref{eq:twomom} are linearly independent on every nonempty interval, and the two source moment functionals are independent.

If both moments vanish, the mean is zero on $O$. It satisfies $Au=0$ on the component of $\mathbb T_L\setminus\supp J$ containing $\mathbb T_L\setminus\overline I$. ODE uniqueness propagates the zero solution through this component and into neighborhoods of both endpoints of $I$. Elliptic regularity gives $u\in C_c^\infty(I)$ and $J=Au$. The converse follows either from \eqref{eq:twomom} by integration by parts or from $A^{-1}Av=v$. A displacement supported in $I$ commutes with every field and momentum Weyl operator supported in $O$, proving equality on the full common time-zero local algebra.
\end{proof}
The proposition holds for any positive mass and the connected separated interval geometry in \eqref{eq:geom}.

For completeness, in physical source coordinates the free global identities are
\begin{equation}\label{eq:BKMfree}
 g_0(f,h)=\beta\ip f{A^{-1}h},\qquad
 D(\rho_{0,J}\Vert\rho_{0,K})=\frac\beta2\norm{J-K}_{H_m^{-1}}^2.
\end{equation}
They follow from $\log Z_0(J)-\log Z_0(0)=\beta\ip J{A^{-1}J}/2$.

\subsection{Thermal covariance and Wick conventions}
Put $k_n=2\pi n/L$, $\omega_n=(k_n^2+m^2)^{1/2}$ and $\nu_j=2\pi j/\beta$. The covariance is
\begin{equation}\label{eq:C}
 C(\tau,x)=\frac1{\beta L}\sum_{j,n\in\Z}
 \frac{e^{i\nu_j\tau+ik_nx}}{\nu_j^2+\omega_n^2}.
\end{equation}
Its time-zero kernel is
\begin{equation}\label{eq:Cbeta}
 C_\beta(x)=\frac1L\sum_{n\in\Z}\frac1{2\omega_n}
 \coth\frac{\beta\omega_n}{2}\,e^{ik_nx}.
\end{equation}
Although $C_\beta(0)$ diverges, $C_\beta\in L^p(\mathbb T_L)$ for every finite $p$. Let $P_N$ truncate only spatial Fourier modes, $\phi_N=P_N\phi$, and $c_N=\E\phi_N(0,x)^2$. Thermal Wick powers use $c_N$ at this cutoff and are passed to the limit as distributions. If $c_{0,N}=L^{-1}\sum_{|n|\le N}(2\omega_n)^{-1}$, then
\begin{equation}\label{eq:dbeta}
 d_\beta=\lim_N(c_N-c_{0,N})
 =\frac1L\sum_{n\in\Z}\frac1{\omega_n(e^{\beta\omega_n}-1)}>0,
\end{equation}
and
\begin{equation}\label{eq:wickreorder}
 :\phi^4:_0=: \phi^4:_\beta+6d_\beta:\phi^2:_\beta+3d_\beta^2.
\end{equation}
Thermal normal ordering defines the model at the fixed inverse temperature $\beta$.

\subsection{Source calibration and the exact shifted measure}
Write $\mathscr T=\mathbb T_\beta\times\mathbb T_L$ and set
\[
 V_4=\int_{\mathscr T}:\phi^4:_\beta,\qquad V_2=\int_{\mathscr T}:\phi^2:_\beta,
 \qquad Y_J=\int_0^\beta\phi(\tau,J)\dd\tau.
\]
The physical Euclidean source measure is proportional to
\begin{equation}\label{eq:physicaldensity}
 \exp\{-\lambda(V_4+bV_2)+Y_{J_{\lambda,v}}\}\dd\mu,
 \qquad b\ge0.
\end{equation}
The associated formal Hamiltonian is $\mathrm d\Gamma(\sqrt A)+\lambda\int(:\Phi^4:_\beta+b:\Phi^2:_\beta)-\Phi(J_{\lambda,v})$; its rigorous KMS interpretation is given in Section~\ref{sec:construct}.

Extend $v\in H=H_0^1(I)$ by zero to the circle. The static Cameron--Martin shift by $v$ has squared norm $\beta\norm v_H^2$. For every local time-zero functional $F$, $F(\phi+v)=F(\phi)$. Completing the quadratic tilt in \eqref{eq:physicaldensity} therefore gives the exact identity
\begin{equation}\label{eq:shift}
 \E_{\lambda,v}F
 =\frac{\E[F e^{-\lambda U(v)}]}{\E[e^{-\lambda U(v)}]},
 \qquad U(v)=\int_{\mathscr T}:(\phi+v)^4:_\beta+bV_2.
\end{equation}
The terms $2b\lambda\int v\phi$ cancel after shifting. The remaining scalar factor is
\[
 \exp\!\left\{\frac\beta2\norm v_H^2+b\lambda\beta\norm v_2^2\right\},
\]
and also cancels in the quotient. We use $\widehat\nu_{\lambda,v}$ for the probability measure on the right of \eqref{eq:shift} and $P_{\lambda,v}^O$ for its local time-zero marginal.

The Gaussian model obtained by deleting $\lambda V_4$ has quadratic operator $A+2b\lambda$, mean $v$, and source-independent covariance. Its local state is thus exactly independent of $v$. The calibration removes the following possible false positive: with the uncalibrated source $Av$, the local mean at infinitesimal source amplitude has the term
\[
 -2b\lambda(A^{-1}v)|_O+O(\lambda^2).
\]
For vacuum ordering $b=6d_\beta$, this is precisely the thermal mass-shift term. This motivates the calibrated source used throughout the interacting model.

The physical source distance and the displacement distance are uniformly comparable:
\begin{equation}\label{eq:sourcecompare}
 \norm{v-w}_H\le\norm{J_{\lambda,v}-J_{\lambda,w}}_{H_m^{-1}}
 \le\left(1+\frac{2b\lambda}{m^2}\right)\norm{v-w}_H.
\end{equation}
This follows by applying the Fourier multiplier $I+2b\lambda A^{-1}$ to $A^{1/2}(v-w)$.

\section{Continuum estimates and the thermal-state interface}\label{sec:construct}
We prove the estimates needed below instead of assuming source-uniform differentiability. The only reconstruction input is the standard positive-temperature FKN theorem, stated at the point where it is used.

\subsection{Wick convergence and Nelson bounds}
Define
\[
 U_N(v)=\int_{\mathscr T}\bigl(:(\phi_N+P_Nv)^4:_{c_N}+b:\phi_N^2:_{c_N}\bigr).
\]
The full Gaussian path space is retained. Only the interaction is truncated; the limiting local algebra is never identified with a spatial Fourier tensor factor.

\begin{lemma}[Uniform polynomial and exponential bounds]\label{lem:Nelson}
For every $R<\infty$, $1\le p<\infty$ and $a_*<\infty$,
\begin{align}
 &\sup_{N,\,\norm v_H\le R}\norm{D_v^j U_N(v)[u_1,\ldots,u_j]}_{L^p(\mu)}
 \le C_{p,R,j}\prod_{i=1}^j\norm{u_i}_H,
 &&0\le j\le4,\label{eq:poly}\\
 &\sup_{N,\,\norm v_H\le R,\,0\le a\le a_*}
 \E e^{-aU_N(v)}<\infty.\label{eq:exp}
\end{align}
The polynomial quantities in \eqref{eq:poly} converge in every finite $L^p$, uniformly on the indicated source and direction balls. Their limits define $U(v)$ and its derivatives. The normalized densities
\begin{equation}\label{eq:rhat}
 \widehat r_{\lambda,v}=\frac{e^{-\lambda U(v)}}{Z_\lambda(v)},
 \qquad Z_\lambda(v)=\E e^{-\lambda U(v)},
\end{equation}
satisfy $\sup_{0\le\lambda\le\lambda_*,\norm v_H\le R}\norm{\widehat r_{\lambda,v}}_2=K_R<\infty$ for every fixed $\lambda_*<\infty$.
\end{lemma}
\begin{proof}
We first establish \eqref{eq:exp} at $v=0$. Let $C_N=P_NC$. The Fourier coefficients in \eqref{eq:C} give
\[
 \norm{C_M-C_N}_{L^2(\mathscr T)}\le CN^{-1},\qquad
 \sup_N\norm{C_N}_{L^6(\mathscr T)}<\infty\quad(M>N).
\]
For the first estimate, sum $(\nu_j^2+\omega_n^2)^{-2}$ over $j$ and $|n|>N$. The second follows from the two-dimensional Sobolev embedding, since $C\in H^a$ for every $a<1$. Wick's rule and nested projections imply
\[
 \left\|\int_{\mathscr T}(:\phi_M^4:-:\phi_N^4:)\right\|_2^2
 =4!\,\beta L\int_{\mathscr T}(C_M^4-C_N^4)
 \le C\norm{C_M-C_N}_2\sup_K\norm{C_K}_6^3\le CN^{-1}.
\]
The same argument for the quadratic term is simpler. Gaussian hypercontractivity on each chaos consequently yields
\begin{equation}\label{eq:tailU}
 \norm{U_M(0)-U_N(0)}_p\le C_b(p-1)^2N^{-1/2},\qquad p\ge2.
\end{equation}

Since $c_N\le C\log(eN)$, the polynomial
\[
 z^4+(b-6c_N)z^2+3c_N^2-bc_N
\]
has a deterministic lower bound $-C_b(1+c_N)^2$. Hence $U_N(0)\ge-C_b\log^2(eN)$. For a sufficiently large $t$, choose $n(t)=\lfloor\exp(\epsilon\sqrt t)\rfloor$ with $\epsilon>0$ small enough that $U_n(0)\ge-t/2$ for $n\le n(t)$. If $N>n(t)$, Markov's inequality and \eqref{eq:tailU} give
\[
 \Pr(U_N(0)<-t)\le
 \left(\frac{2C_b(p-1)^2n(t)^{-1/2}}t\right)^p.
\]
Taking an integer $p$ proportional to $\sqrt t\,n(t)^{1/4}$, with a sufficiently small proportionality constant, proves
\[
 \sup_N\Pr(U_N(0)<-t)\le
 \exp\{-c\sqrt t\,e^{c'\sqrt t}\}.
\]
This proves every negative exponential moment uniformly in $N$. Fatou's lemma includes the limit, and
\[
 |e^{-ax}-e^{-ay}|\le a|x-y|(e^{-ax}+e^{-ay})
\]
with H\"older's inequality proves convergence of weights in every finite $L^p$, uniformly for bounded $a\ge0$. This is the classical Nelson argument, included to specify the cutoff and constants' dependencies; see \cite{Simon,GlimmJaffe} for its general constructive context.

For source dependence, $H^1(\mathbb T_L)\hookrightarrow L^\infty$ and $\norm{P_Nv}_\infty\le C\norm v_H$. Expansion of $:(\phi_N+P_Nv)^4:$ gives finitely many Wick polynomials with bounded deterministic coefficients. If $j\le4$ and $g$ is bounded, Wick's rule gives
\[
 \norm{\int_{\mathscr T}g:\phi_N^j:}_2^2
 \le j!\norm g_\infty^2\beta L\norm{C_N^j}_{L^1}.
\]
Hypercontractivity proves \eqref{eq:poly}. For convergence use $C_N^j\to C^j$ in $L^1$, together with
$\norm{(I-P_N)v}_\infty\le CN^{-1/2}\norm v_H$.

To obtain the shifted exponential bound, apply the static Cameron--Martin identity with the \emph{complete} shift $v$. It gives
\begin{equation}\label{eq:CMbound}
 \E e^{-aU_N(v)}=
 e^{-\beta\norm v_H^2/2-ab\beta\norm{P_Nv}_2^2}
 \E e^{-aU_N(0)+Y_{Av+2abP_Nv}}.
\end{equation}
The prefactor is at most one. Cauchy--Schwarz bounds the last expectation by
\[
 (\E e^{-2aU_N(0)})^{1/2}
 \exp\!\left\{\beta\left(1+\frac{2ab}{m^2}\right)^2R^2\right\}.
\]
This proves \eqref{eq:exp}. Finally $\E U_N(v)=\beta\int(P_Nv)^4\le C R^4$, so Jensen's inequality gives $Z_{N,\lambda}(v)\ge e^{-C\lambda R^4}$. Dividing the uniform exponential bounds by this lower bound proves the density estimate and its convergence.
\end{proof}

\subsection{Local cubic fields and the quantum realization}
\begin{lemma}[Sharp-time observables]\label{lem:Wh}
For $h\in L^2(O)$ the thermal Wick cubic $W_h$ exists in every finite $L^p(\mu)$ and is measurable with respect to $\OO=\sigma\{\phi(0,h):h\in C_c^\infty(O)\}$. Moreover
\begin{align}
 \E W_hW_k&=6\int_{O\times O}h(y)k(z)C_\beta(y-z)^3\dd y\dd z,\label{eq:Wcov}\\
 \norm{W_h}_p&\le c_p\norm h_2,\qquad
 \sup_{\lambda\le\lambda_*,\norm v_H\le R}\E_{\widehat\nu_{\lambda,v}}W_h^2
 \le C_W\norm h_2^2.\label{eq:Wmom}
\end{align}
One admissible choice is $C_W=162K_R\norm{C_\beta^3}_{L^1(\mathbb T_L)}$.
\end{lemma}
\begin{proof}
The Fourier coefficients of $C_\beta$ have size $O((1+|n|)^{-1})$, hence $C_\beta\in L^q$ for all finite $q$, and the projected covariances converge in each such space. The covariance of the difference of two nested Wick cubics is $6(C_{\beta,M}^3-C_{\beta,N}^3)$, which tends to zero in $L^1$. Young's inequality gives convergence as maps $h\mapsto W_h$ from $L^2$ to Gaussian $L^2$. Hypercontractivity upgrades this to every finite $L^p$ and proves \eqref{eq:Wcov}.

For locality first take $h\in C_c^\infty(O)$ and use local mollifiers with support smaller than the distance of $\supp h$ to $\partial O$. The same covariance argument shows that their Wick cubics have the same limit as the Fourier construction. They are $\OO$-measurable. Density then treats all $h\in L^2(O)$. 

Since $W_h$ lies in third Gaussian chaos, $\norm{W_h}_4\le3^{3/2}\norm{W_h}_2$. Consequently
\[
 \E\widehat r_{\lambda,v}W_h^2\le K_R\norm{W_h}_4^2
 \le27K_R\E W_h^2\le162K_R\norm{C_\beta^3}_1\norm h_2^2.
\]
\end{proof}

\begin{proposition}[Thermal-state interface]\label{prop:FKN}
For $v\in C_c^\infty(I)$, the finite-volume free thermal system admits the FKN perturbation by the sharp-time potential
\[
 \mathscr V_{\lambda,v}=\lambda\int_{\mathbb T_L}(:\Phi^4:_\beta+b:\Phi^2:_\beta)\dd x-\Phi(J_{\lambda,v}).
\]
The resulting states $\omega_{\lambda,v}$ are normal states on a common observable algebra. Their restrictions to its time-zero field subalgebra have the periodic Euclidean law \eqref{eq:physicaldensity}. In particular their equal-time restrictions to $O$ have law $P_{\lambda,v}^O$ from \eqref{eq:shift}.
\end{proposition}
\begin{proof}
We first identify precisely the perturbation theorem. Theorem~7.11 of
G\'erard--J\"akel \cite{GJcut} applies to a stochastically positive
$\beta$-KMS system $(\cF,\mathcal U,\tau,\omega)$ and to a real,
selfadjoint time-zero potential $V$ affiliated with its abelian algebra
$\mathcal U$, provided that
\begin{equation}\label{eq:FKNassumptions}
 V,e^{-\beta V}\in L^1(Q,\Sigma_0,\mu),\qquad
 V\in L^{2+\epsilon}(Q,\Sigma_0,\mu)\quad\text{for some }\epsilon>0.
\end{equation}
We apply the $L^{2+\epsilon}$ integrability criterion. Theorem~7.6 establishes the periodic
FKN path measure, while Propositions~7.8--7.10 identify the time-zero
abelian algebra, the observable algebra, and the perturbed dynamics under
the canonical unitary identification.

For clarity, we verify the free compact-space input. With
$k_n=2\pi n/L$ and $\omega_n=(k_n^2+m^2)^{1/2}$, its time-periodic
covariance is, for $|t|\le\beta$,
\begin{equation}\label{eq:OScircle}
 C(t,x)=\frac1L\sum_{n\in\Z}e^{ik_nx}
 \frac{e^{-\omega_n|t|}+e^{-\omega_n(\beta-|t|)}}
 {2\omega_n(1-e^{-\beta\omega_n})}.
\end{equation}
Each spatial Fourier mode is the covariance of a thermal harmonic
oscillator. In particular, for $0\le s,t\le\beta/2$, the time-reflected
kernel of one mode is
\[
 \frac{e^{-\omega_n(s+t)}+e^{-\beta\omega_n}e^{\omega_n(s+t)}}
 {2\omega_n(1-e^{-\beta\omega_n})},
\]
a sum of two positive rank-one kernels. Summing the modes gives
reflection positivity for the spatially smeared fields. The Gaussian
bridge conditional-independence property for each oscillator, followed
by the finite-mode limit, gives the periodic Markov property and the
standard stochastically positive free KMS path space. Equivalently, the
circle heat semigroup is the explicitly defined spectral semigroup
\[
 (e^{t\partial_x^2}f)(x)=\frac1L\sum_{n\in\Z}
    e^{-tk_n^2}\widehat f(n)e^{ik_nx},\qquad t>0,
\]
which requires no noncompact-space estimate; see \cite{Davies,GJcut,KleinLandau}.
The general results of Section~7 of \cite{GJcut} are formulated for the
abstract KMS path space and do not require the spatial coordinate to be
$\R$ once these free-system hypotheses are checked.

Apply the theorem with
\[
 V=\mathscr V_{\lambda,v}=\lambda\int_{\mathbb T_L}
       \bigl(:\Phi^4:_\beta+b:\Phi^2:_\beta\bigr)\dd x
       -\Phi(J_{\lambda,v}),\qquad v\in C_c^\infty(I).
\]
At sharp time, $\norm{C_{\beta,M}-C_{\beta,N}}_2\le CN^{-1/2}$ and
$\sup_N\norm{C_{\beta,N}}_6<\infty$. Repeating the Wick-contraction
and hypercontractivity argument in Lemma~\ref{lem:Nelson}, with a single
spatial slice instead of the full time--space cylinder, gives
$L^p$ convergence of the integrated Wick quartic for all finite $p$,
with error $C_pN^{-1/4}$. The deterministic cutoff lower bound is
$-C\log^2(eN)$, and the same low/high-frequency Nelson argument gives
\begin{equation}\label{eq:sharptimeexp}
 \sup_N\E\exp\bigl[-q\beta\lambda
    \textstyle\int_{\mathbb T_L}(:\Phi_N^4:_\beta
       +b:\Phi_N^2:_\beta)\dd x\bigr]<\infty,
 \quad 0\le q<\infty,\ \lambda\ge0.
\end{equation}
The smooth source $\Phi(J_{\lambda,v})$ is a centered Gaussian linear
functional at sharp time, hence $\E e^{q\Phi(J_{\lambda,v})}<\infty$
for every real $q$. H\"older's inequality, the $L^p$ convergence of
the Wick integral, and \eqref{eq:sharptimeexp} imply that
$V\in L^p$ for every finite $p$ and $e^{-\beta V}\in L^1$.
Both the Wick integral and the smooth field functional are measurable
with respect to the sharp-time sigma algebra; their limits therefore
define a real affiliated multiplication operator. This verifies each
hypothesis in \eqref{eq:FKNassumptions} for the circle model.

The normalized Euclidean path measure supplied by the periodic FKN
theorem has density \eqref{eq:physicaldensity}. Its sharp-time law,
translated by the static Cameron--Martin shift, is exactly
\eqref{eq:shift}. The common normal-state realization and the stated
equal-time local restrictions follow from Theorem~7.11 and
Propositions~7.8--7.10 of \cite{GJcut}.
\end{proof}
We compare all states on a fixed time-zero local algebra $\cA_0(O)$ containing the Weyl field and momentum operators supported in $O$. Write
\[
 S_O(\omega_{\lambda,v}\Vert\omega_{\lambda,w})
 =S(\omega_{\lambda,v}|_{\cA_0(O)}\Vert\omega_{\lambda,w}|_{\cA_0(O)}).
\]
Every $W_h$ is affiliated to its time-zero field subalgebra. Monotonicity of Araki relative entropy \cite{Araki} gives
\begin{equation}\label{eq:dataproc}
 S_O(\omega_{\lambda,v}\Vert\omega_{\lambda,w})\ge
 \KL(P_{\lambda,v}^O\Vert P_{\lambda,w}^O).
\end{equation}
For nonsmooth $v\in H$, the probability estimates refer to \eqref{eq:shift}; the quantum interpretation applies to smooth sources as in Proposition~\ref{prop:FKN}.

\section{A uniform cubic-response expansion}\label{sec:expansion}
\begin{theorem}[Source-uniform weak-coupling response]\label{thm:response}
Let $b\ge0$ and $\lambda_*<\infty$ be fixed. For every $R<\infty$ the functionals
\[
 h\longmapsto \E_{\widehat\nu_{\lambda,v}}W_h
\]
define $\mathcal M_\lambda(v)\in Y=L^2(O)$ for $v\in H$, $\norm v_H\le R$, and $0\le\lambda\le\lambda_*$. There is a constant $C_R$ such that
\begin{equation}\label{eq:expansion}
 \mathcal M_\lambda(v)=-24\lambda\T v+\mathcal R_\lambda(v),\qquad
 \mathcal R_\lambda(0)=0,\qquad
 \sup_{\norm v_H\le R}\norm{D_v\mathcal R_\lambda(v)}_{H\to Y}\le C_R\lambda^2.
\end{equation}
All statements are about the continuum model. At finite spatial cutoff the analogous estimate holds with its own leading kernel $\T_N$, with constants independent of the cutoff. Passing to the limit gives \eqref{eq:expansion} with the full kernel $\T$.
\end{theorem}
\begin{proof}
Put $F_h(\lambda,v)=\E[W_he^{-\lambda U(v)}]$ and $Z(\lambda,v)=\E e^{-\lambda U(v)}$. Riesz representation and Lemmas~\ref{lem:Nelson}--\ref{lem:Wh} define $\mathcal M_\lambda$ because
\[
 |F_h/Z|\le K_R\norm{W_h}_2\le C_R\norm h_2.
\]
We justify, in particular, the mixed derivative required for the remainder, Let $Q_v[u]=D_vU(v)[u]$. For $j\le2$,
\[
 \partial_\lambda^jF_h=\E[W_h(-U(v))^je^{-\lambda U(v)}],\qquad
 D_vF_h[u]=-\lambda\E[W_hQ_v[u]e^{-\lambda U(v)}].
\]
In particular,
\begin{equation}\label{eq:mixednumerator}
 \partial_\lambda^2D_vF_h[u]
 =\E[W_hQ_v[u](2U(v)-\lambda U(v)^2)e^{-\lambda U(v)}].
\end{equation}
The same formulas with $W_h$ omitted hold for $Z$. For every fixed product in these expressions, H\"older's inequality, \eqref{eq:poly}, \eqref{eq:exp}, and \eqref{eq:Wmom} give a uniform bound $C_R\norm h_2\norm u_H$. Taylor's formula for the polynomial $U(v+u)$ proves Fr\'echet differentiability. The same estimates dominate the difference quotients. Since $Z\ge e^{-C\lambda_*R^4}$, the quotient rule gives
\begin{equation}\label{eq:mixedbound}
 \sup_{\lambda\in[0,\lambda_*],\norm v_H\le R}
 |\partial_\lambda^2D_v(F_h/Z)[u]|
 \le C_R\norm h_2\norm u_H.
\end{equation}
The cutoff convergence in Lemma~\ref{lem:Nelson}, again with H\"older, passes all these derivatives to the continuum uniformly on the balls. The derivatives in $\lambda$ at zero are right derivatives; no convergent perturbation series is required.

At $\lambda=0$, $\mathcal M_0(v)=0$. Expanding the shifted Wick polynomial gives
\begin{equation}\label{eq:Uexpand}
 U(v)=V_4+bV_2+4\int_{\mathscr T}v:\phi^3:
 +6\int_{\mathscr T}v^2:\phi^2:
 +4\int_{\mathscr T}v^3\phi+\beta\int_Iv^4.
\end{equation}
Only the third-chaos term pairs with $W_h$. Wick's rule yields
\begin{align*}
 \partial_\lambda\ip h{\mathcal M_\lambda(v)}\big|_{0+}
 &=-\E[W_hU(v)]\\
 &=-24\int_{O\times I}h(y)v(x)
       \left(\int_0^\beta C(\tau,x-y)^3\dd\tau\right)\dd x\dd y
 =-24\ip h{\T v}.
\end{align*}
Taylor's integral formula applied to $D_v\mathcal M_\lambda$ and \eqref{eq:mixedbound} proves \eqref{eq:expansion}. The measure at $v=0$ is invariant under $\phi\mapsto-\phi$, so its cubic mean is exactly zero and $\mathcal R_\lambda(0)=0$.
\end{proof}

\begin{remark}[What the cubic observable removes]\label{rem:calibration}
The first coefficient in \eqref{eq:expansion} is unchanged by $bV_2$ or by an additional first-order quadratic counterterm accompanied by the corresponding quadratic-reference source calibration. Such terms have Gaussian chaos degree at most two and are orthogonal to $W_h$. Thus the coefficient is not a thermal mass-shift artifact. The finite-coupling expectation $\E W_h$ is a Wick cubic moment, not in general the exact connected third cumulant. Their first-order non-Gaussian response agrees here because the local free mean is zero. 
\end{remark}

For comparison, for a local linear test $h$ the same calculation gives, with $v=su$,
\begin{equation}\label{eq:calmean}
 \E_{\lambda,su}\Phi(h)
 =-4\lambda s^3\int_Iu(x)^3(A^{-1}h)(x)\dd x+O_{C^1([-S,S])}(\lambda^2).
\end{equation}
The term linear in $s$ has disappeared after the calibration. Its derivative at $s=0$ therefore cannot furnish the uniform first-order lower gain that an uncalibrated thermal mass term would supply. The cubic channel \eqref{eq:expansion}, by contrast, is linear in the entire displacement $v$ at its first nonzero order.

\section{The thermal cubic kernel: injectivity, rank, and ill-posedness}\label{sec:kernel}
\subsection{A positive bilateral Laplace representation}
Set $\kappa_n=(m^2+(2\pi n/\beta)^2)^{1/2}$. On the cut circle $0<z<L$, the Matsubara expansion and its time integral give
\begin{align}
 C(\tau,z)&=\frac1\beta\sum_{n\in\Z}e^{2\pi in\tau/\beta}G_{\kappa_n}(z),\label{eq:mats}\\
 \K(z)&=\frac1{\beta^2}\sum_{n,j\in\Z}
 G_{\kappa_n}(z)G_{\kappa_j}(z)G_{\kappa_{n+j}}(z).\label{eq:Ksum}
\end{align}
For every closed substrip $\epsilon\le\Re z\le L-\epsilon$, the double sum and all its fixed-order derivatives converge absolutely and uniformly. Indeed each factor and its derivatives are bounded by a polynomial in $\kappa$ times $e^{-\epsilon\kappa}/\kappa$.

\begin{lemma}[Positive atoms]\label{lem:atoms}
There exists a positive, locally finite atomic measure $\sigma_\beta$ on $\R$ such that
\begin{equation}\label{eq:Laplace}
 \K(z)=\int_\R e^{-az}\dd\sigma_\beta(a),\qquad 0<\Re z<L.
\end{equation}
It has an atom of mass at least
\begin{equation}\label{eq:atoms}
 b_n=\frac1{8\beta^2m\kappa_n^2(1-e^{-mL})(1-e^{-\kappa_nL})^2}
\end{equation}
at $a_n=2\kappa_n+m$, for every $n\ge1$. In particular $b_n\ge c(1+n)^{-2}$, $a_n\asymp n$, and $\inf_n(a_{n+1}-a_n)>0$.
\end{lemma}
\begin{proof}
Expand each Green function in \eqref{eq:Ksum} into its two positive exponential terms. This defines a positive atomic measure, allowing repeated exponents to add their masses. Finiteness of \eqref{eq:Ksum} at any real $z\in(0,L)$ proves local finiteness and the bilateral exponential integrability required in \eqref{eq:Laplace}. The term $(n,j)=(n,0)$, with the decaying exponential chosen in all three factors, gives \eqref{eq:atoms}. The remaining assertions follow directly from the formula for $\kappa_n$.
\end{proof}

\begin{theorem}[Infinite-dimensional injectivity]\label{thm:injective}
For the separated geometry \eqref{eq:geom}, $\T:L^2(I)\to L^2(O)$ is injective.
\end{theorem}
\begin{proof}
Suppose $\T u=0$ on $O$. For $\ell-L<\Re y<0$ define
\[
 F(y)=\int_I\K(x-y)u(x)\dd x,
 \qquad U(a)=\int_Ie^{-ax}u(x)\dd x.
\]
The kernel estimates above show that $F$ is holomorphic in this vertical strip. Its continuous restriction vanishes on the real interval $O$, so $F$ vanishes on the entire strip. Fix a real $y_0$ in the strip. The signed measure
\[
 \dd\zeta(a)=e^{ay_0}U(a)\dd\sigma_\beta(a)
\]
has finite total variation: its integral is bounded by
$\int_I|u(x)|\K(x-y_0)\dd x<\infty$. For real $t$,
\[
 F(y_0+it)=\int_\R e^{iat}\dd\zeta(a)=0.
\]
Fourier uniqueness for finite measures gives $\zeta=0$. Positivity of the atom at each $a_n$ implies $U(a_n)=0$. 

The function $U(z)=\int_0^\ell e^{-zx}u(x)\dd x$ is bounded and holomorphic on the right half-plane. If it is nonzero, its positive real zeros obey the Blaschke condition
$\sum a_n/(1+a_n^2)<\infty$. To recall the elementary reason, map the half-plane to the disk by $\zeta=(z-1)/(z+1)$ and apply Jensen's formula to the resulting bounded holomorphic function (after removing a possible finite-order zero at the center). Its zeros must have $\sum(1-|\zeta_n|)<\infty$, whereas for $a_n\asymp n$ this sum diverges. Therefore $U$ is identically zero. Fourier uniqueness for $e^{-cx}u(x)$ with any $c>0$ gives $u=0$. Equivalently, the last step is the M\"untz completeness mechanism \cite{BorweinErdelyi}; the half-plane proof specifies the zero-set condition needed here.
\end{proof}

\subsection{Finite observable designs and an analytic rank distinction}
\begin{proposition}[Strict sign regularity]\label{prop:sign}
If $y_1<\cdots<y_r<x_1<\cdots<x_r$ and $x_r-y_1<L$, then
\begin{equation}\label{eq:sign}
 (-1)^{r(r-1)/2}\det[\K(x_j-y_i)]_{i,j=1}^r>0.
\end{equation}
Consequently nonzero nonnegative smooth probes and source displacement profiles supported on respectively ordered, disjoint subintervals give a full-rank cubic response matrix of every fixed finite size.
\end{proposition}
\begin{proof}
The continuous determinant identity applied to \eqref{eq:Laplace} gives
\[
 \det[\K(x_j-y_i)]
 =\frac1{r!}\int_{\R^r}
 \det[e^{a_ky_i}]_{i,k}\det[e^{-a_kx_j}]_{k,j}
 \prod_{k=1}^r\dd\sigma_\beta(a_k).
\]
Absolute integrability follows by expanding the determinants and using the finite positive kernel products. For ordered distinct $a_k$, the first determinant is positive and the second has sign $(-1)^{r(r-1)/2}$. This sign property of exponentials follows from their nonzero Wronskian and Rolle's theorem; its sign is fixed by coalescing the ordered arguments. The measure has infinitely many distinct atoms, so a set with $r$ distinct atoms has positive measure. This proves strictness. Integrating \eqref{eq:sign} against the ordered nonnegative profiles proves the last assertion by multilinearity.
\end{proof}

For uncalibrated linear field means, the first-order mass-shift response on one exterior interval lies in the span of $e^{my}$ and $e^{-my}$, so it has rank at most two. If only the zero Matsubara mode is kept in \eqref{eq:Ksum}, the cubic kernel is $\beta^{-2}G_m(z)^3$ and its rank is at most four, from the exponents $\pm m,\pm3m$. The full thermal kernel has infinite rank by Theorem~\ref{thm:injective}. Effective rank additionally depends on the observation noise, spatial smoothing and sample size.

\begin{proposition}[Analytic smoothing]\label{prop:compact}
The operators $\T:L^2(I)\to L^2(O)$ and $\T:H\to L^2(O)$ are compact. There are $C,c>0$ such that their singular values satisfy $s_n(\T)\le Ce^{-cn}$. No bound $\norm u_H\le C\norm{\T u}_Y$ holds on $H$. On any nontrivial spectral prior ball $\norm u_{X_s}\le M$, no uniform H\"older inverse estimate $\norm u_H\le C\norm{\T u}_Y^\alpha$ with $\alpha>0$ holds.
\end{proposition}
\begin{proof}
Because $\overline I$ and $\overline O$ are separated, $\K(x-y)$ is holomorphic and bounded in a fixed complex neighborhood of $y\in\overline O$, uniformly in $x\in\overline I$. Rescale the observation interval to $[-1,1]$ and choose a Bernstein ellipse of parameter $\rho>1$ lying in that common neighborhood. The Chebyshev coefficients $c_j(x)$ of $y\mapsto\K(x-y)$ satisfy
\[
 \sup_{x\in I}|c_j(x)|\le C\rho^{-j}.
\]
This follows by writing the coefficient as a contour integral on the ellipse (the classical analytic approximation estimate \cite[Chap.~8]{Trefethen}). Thus the truncated kernel
$\sum_{j=0}^{n-1}c_j(x)T_j(\widetilde y)$ is a sum of $n$ separated functions of $x$ and $y$, hence gives an operator of rank at most $n$. The uniformly convergent geometric tail has supremum norm at most $C\rho^{-n}$ on $I\times O$ and therefore Hilbert--Schmidt and operator norm at most $C'\rho^{-n}$. The approximation-number characterization of Hilbert-space singular values gives $s_{n+1}(\T)\le C'\rho^{-n}$. Finally, composition with the bounded embedding $H\hookrightarrow L^2(I)$ gives the same form for the map from $H$ to $Y$.

For the last assertion, let $E_n$ be the first $n$ Dirichlet modes, normalized in $H$. A rank-$(n-1)$ approximant has a unit vector $u_n\in E_n$ in its kernel, whence $\norm{\T u_n}\le Ce^{-cn}$ and $\norm{u_n}_{X_s}\le Cn^s$. Scaling to the prior ball produces $\norm{v_n}_H\ge cMn^{-s}$ and $\norm{\T v_n}\le CMn^{-s}e^{-cn}$. These contradict any such H\"older estimate. Compactness and infinite-dimensional injectivity also rule out a bounded inverse on the whole space.
\end{proof}

\subsection{A reproducible finite-grid singular-value illustration}
\begin{figure}[!htbp]
\centering
\includegraphics[width=.68\linewidth]{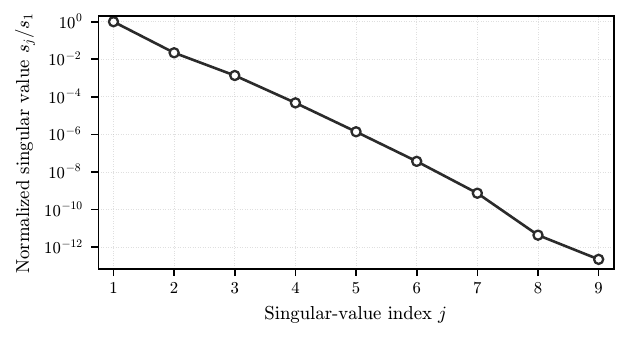}
\caption{Normalized singular values of the $96\times96$ weighted finite-grid matrix \eqref{eq:nummatrix}. The thermal Matsubara sum is truncated at $|n|\le35$; the independent $48\times48$ comparison is reported in the text. Reproduction code and numerical data are included in the accompanying ancillary files.}
\label{fig:singvalues}
\end{figure}

Take $\beta=4$, $m=1$, $L=2\pi$, $I=(0,1)$ and $O=(-1.8,-1.0)$. These intervals satisfy \eqref{eq:geom} and the short-source condition \eqref{eq:short}. Using the Matsubara expansion \eqref{eq:mats}, retain $|n|\le35$ and approximate the periodic $\tau$-integral in $\K(z)=\int_0^\beta C(\tau,z)^3\dd\tau$ by a 192-point midpoint trapezoidal rule. On uniform midpoint grids $x_j\in I$ and $y_i\in O$ of size $N$, form the weighted matrix
\begin{equation}\label{eq:nummatrix}
 A^{(N)}_{ij}=\sqrt{\Delta x\,\Delta y}\,\K(x_j-y_i).
\end{equation}
Figure~\ref{fig:singvalues} displays the normalized singular values for $N=96$. Their values at indices $2,4,6,8$ are approximately $2.22\cdot10^{-2}$, $4.79\cdot10^{-5}$, $3.69\cdot10^{-8}$, and $4.35\cdot10^{-12}$. An independent $N=48$ calculation differs by at most about $2.2\%$ for the first six normalized singular values. These finite-matrix singular values illustrate Proposition~\ref{prop:compact}; they are not certified bounds for the continuum operator and the smallest values remain sensitive to floating-point precision.
\FloatBarrier
\Needspace{6\baselineskip}

\section{Local entropy bounds and Sobolev resolution}\label{sec:entropy}
We use the convention $H^2(P,Q)=\int(\sqrt{\dd P}-\sqrt{\dd Q})^2$. It satisfies $\KL(P\Vert Q)\ge H^2(P,Q)$.

\Needspace{9\baselineskip}
\begin{lemma}[A second-moment information bound]\label{lem:info}
For the local laws defined above and $\norm v_H,\norm w_H\le R$,
\begin{equation}\label{eq:info}
 \KL(P_{\lambda,v}^O\Vert P_{\lambda,w}^O)
 \ge\frac{\norm{\mathcal M_\lambda(v)-\mathcal M_\lambda(w)}_Y^2}{4C_W}.
\end{equation}
For smooth sources the same right side is a lower bound for $S_O(\omega_{\lambda,v}\Vert\omega_{\lambda,w})$.
\end{lemma}
\begin{proof}
For any $h$, factor the difference of densities as $(\sqrt p-\sqrt q)(\sqrt p+\sqrt q)$. Cauchy--Schwarz and \eqref{eq:Wmom} yield
\[
 |\E_PW_h-\E_QW_h|^2
 \le 2(\E_PW_h^2+\E_QW_h^2)H^2(P,Q)
 \le4C_W\norm h_2^2H^2(P,Q).
\]
Taking the supremum over unit $h$ proves \eqref{eq:info}; smooth compactly supported tests are dense. The quantum statement is \eqref{eq:dataproc}. 
\end{proof}

\begin{theorem}[Fixed finite-dimensional source spaces]\label{thm:finite}
Let $E\subset H$ be any fixed finite-dimensional subspace and set
\[
 \gamma_E=\inf_{u\in E,\norm u_H=1}\norm{\T u}_Y>0.
\]
For $v,w\in E$ of $H$-norm at most $R$ and
$0<\lambda\le\min\{\lambda_*,12\gamma_E/C_R\}$ (omit the second restriction if $C_R=0$),
\begin{equation}\label{eq:finite}
 \KL(P_{\lambda,v}^O\Vert P_{\lambda,w}^O)
 \ge\frac{36\lambda^2\gamma_E^2}{C_W}\norm{v-w}_H^2.
\end{equation}
If $E\subset C_c^\infty(I)$, this also holds with $S_O$ on the left. By \eqref{eq:sourcecompare} the right side is at least
\[
 \frac{36\lambda^2\gamma_E^2}{C_W(1+2b\lambda/m^2)^2}
 \norm{J_{\lambda,v}-J_{\lambda,w}}_{H_m^{-1}}^2.
\]
\end{theorem}
\begin{proof}
Injectivity and compactness of the unit sphere of $E$ give $\gamma_E>0$. The $C^1$ estimate in Theorem~\ref{thm:response}, integrated on the source segment, yields
\[
 \norm{\mathcal M_\lambda(v)-\mathcal M_\lambda(w)}
 \ge\lambda(24\gamma_E-C_R\lambda)\norm{v-w}_H
 \ge12\lambda\gamma_E\norm{v-w}_H.
\]
Apply Lemma~\ref{lem:info}. 
\end{proof}

\subsection{An infinite-dimensional prior at prescribed resolution}
Let $\mu_k=(\pi k/\ell)^2+m^2$ and
\[
 e_k(x)=(2/\ell)^{1/2}\mu_k^{-1/2}\sin(\pi kx/\ell),\qquad De_k=\mu_ke_k.
\]
These form an orthonormal basis of $H$. For $s>0$ define
\begin{equation}\label{eq:Xs}
 X_s=\mathcal D(D^{s/2}),\qquad
 \norm v_{X_s}=\norm{D^{s/2}v}_H,\qquad
 \mathcal B_{M,s}=\{v\in X_s:\norm v_{X_s}\le M\}.
\end{equation}
For larger $s$, the Dirichlet spectral definition fixes the boundary compatibility conditions. The ball is compact in $H$, and its $H$-radius is at most $R=\mu_1^{-s/2}M$.

\begin{theorem}[Uniform separation at every prescribed resolution]\label{thm:resolution}
Fix $M,s>0$. For $r>0$ such that the set below is nonempty, put
\begin{equation}\label{eq:alpha}
 \alpha_{M,s}(r)=
 \min_{\substack{v,w\in\mathcal B_{M,s}\\\norm{v-w}_H\ge r}}
 \norm{\T(v-w)}_Y.
\end{equation}
Then $\alpha_{M,s}(r)>0$. For every
\begin{equation}\label{eq:lambdar}
 0<\lambda\le\lambda_r:=
 \min\left\{\lambda_*,\frac{6\alpha_{M,s}(r)}{C_RR}\right\}
\end{equation}
(with the second restriction omitted if $C_R=0$), all pairs in this set satisfy
\begin{equation}\label{eq:resolution}
 \KL(P_{\lambda,v}^O\Vert P_{\lambda,w}^O)
 \ge\frac{36\lambda^2\alpha_{M,s}(r)^2}{C_W}>0.
\end{equation}
For its smooth pairs the same bound holds for the local Araki relative entropy.
\end{theorem}
\begin{proof}
The minimizing set in \eqref{eq:alpha} is compact in $H\times H$. A zero minimum would give $v=w$ by Theorem~\ref{thm:injective}, contradicting the resolution constraint. Since $\norm{v-w}_H\le2R$, Theorem~\ref{thm:response} gives
\[
 \norm{\mathcal M_\lambda(v)-\mathcal M_\lambda(w)}
 \ge24\lambda\alpha_{M,s}(r)-2C_RR\lambda^2
 \ge12\lambda\alpha_{M,s}(r).
\]
Lemma~\ref{lem:info} proves \eqref{eq:resolution}.
\end{proof}
Here $\lambda_r$ depends on the prescribed resolution $r$; fixed-coupling separation at all resolutions is Problem~\ref{prob:fixed}.

\subsection{The weak-coupling information coefficient}
\begin{proposition}[Non-Gaussian local information at order $\lambda^2$]\label{prop:Hellinger}
For fixed $v,w\in H$, set
\[
 B_{v,w}=-\E[U(v)-U(w)\mid\OO]+\E[U(v)-U(w)].
\]
Then
\begin{equation}\label{eq:Hellinger}
 \lim_{\lambda\downarrow0}\lambda^{-2}H^2(P_{\lambda,v}^O,P_{\lambda,w}^O)
 =\tfrac14\norm{B_{v,w}}_2^2
 \ge144\sup_{h\ne0}
 \frac{|\ip h{\T(v-w)}|^2}{\E W_h^2}.
\end{equation}
The supremum ranges over nonzero $h\in C_c^\infty(O)$ with nonzero denominator. For $v\ne w$ it is positive. For smooth sources the right side is also a lower bound for $\liminf_{\lambda\downarrow0}S_O/\lambda^2$.
\end{proposition}
\begin{proof}
Relative to the free local law, the density is
$p_{\lambda,v}=\E[e^{-\lambda U(v)}\mid\OO]/Z_\lambda(v)$.
Lemma~\ref{lem:Nelson} gives its expansion in $L^4$:
\[
 p_{\lambda,v}=1-\lambda\{\E[U(v)\mid\OO]-\E U(v)\}+O_{L^4}(\lambda^2).
\]
For $x\ge-1$, the error in $\sqrt{1+x}=1+x/2+O(x^2)$ has a uniform absolute bound $x^2/2$. Thus the square-root densities have the corresponding $L^2$ expansion, proving the equality in \eqref{eq:Hellinger}. Equation~\eqref{eq:Uexpand} and Wick orthogonality give
$\E[W_hB_{v,w}]=-24\ip h{\T(v-w)}$.
Cauchy--Schwarz proves the lower bound; injectivity proves its strict positivity for a suitable test. Finally use \eqref{eq:dataproc} and $\KL\ge H^2$.
\end{proof}
Equation~\eqref{eq:Hellinger} gives the measured Hellinger coefficient and a lower coefficient for local Araki entropy; the corresponding classical Fisher information and quantum BKM metric describe different quantities.

\section{An explicit conditional stability estimate}\label{sec:stability}
Theorem~\ref{thm:resolution} needs only compactness and injectivity. We now obtain an explicit, nonoptimal separation function for the actual kernel under the additional geometric condition
\begin{equation}\label{eq:short}
 0<\ell<L/2.
\end{equation}

\subsection{Reflection energy and analytic propagation}
Choose $d>0$ with $2\ell+d<L$ and set $O^\sharp=(-d-\ell,-d)$. This interval is an analytic-continuation target, not an additional measurement region. If $F=\T u$, define
\begin{equation}\label{eq:reflection}
 \mathcal E_d(u)=\int_Iu(x')F(-d-x')\dd x'
 =\int_\R e^{-ad}\left|\int_Ie^{-ax}u(x)\dd x\right|^2\dd\sigma_\beta(a)\ge0.
\end{equation}
The identity follows from Fubini and \eqref{eq:Laplace}, since $0<d+x+x'<L$. Absolute integrability holds even for sign-changing $u$, by the boundedness of the kernel on this compact separation range. The positivity of the reflected quadratic energy permits restriction to any subset of positive atoms for sign-changing $u$.

\begin{lemma}[Propagation to the reflected interval]\label{lem:propagate}
There are $C>0$ and $\alpha\in(0,1)$, depending only on the fixed geometry and thermal parameters, such that, for $\norm u_H\le B$ and $\delta=\norm{\T u}_{L^2(O)}\le B$,
\begin{equation}\label{eq:propagate}
 \sup_{y\in\overline{O^\sharp}}|F(y)|\le CB^{1-\alpha}\delta^\alpha.
\end{equation}
Consequently, with $U_n=\int_Ie^{-a_nx}u(x)\dd x$,
\begin{equation}\label{eq:moments}
 \sum_{n=1}^\infty b_ne^{-da_n}|U_n|^2\le CB^{2-\alpha}\delta^\alpha.
\end{equation}
\end{lemma}
\begin{proof}
The analyticity statements in Section~\ref{sec:kernel} imply that $F$
is holomorphic throughout the strip
$S=\{z\in\mathbb C:\ell-L<\Re z<0\}$. Since the closures of $O$
and $O^\sharp$ lie strictly inside the real part of $S$, choose a
bounded open rectangle $D_0\Subset S$ whose interior contains both,
and a slightly larger rectangle $D_1\Subset S$ with
$\overline{D_0}\subset D_1$. Uniform kernel bounds on
$\overline{D_1}\times\overline I$ and Cauchy--Schwarz give
\begin{equation}\label{eq:propM}
 \sup_{z\in\overline{D_1}}|F(z)|+
 \sup_{z\in\overline{D_0}}|F'(z)|\le C_0B.
\end{equation}
Fix closed real intervals $J_0\Subset\operatorname{int}J_1\Subset O$
with positive lengths and put $A_0=\max_{J_0}|F|$.
Assume first $B>0$ and $\delta>0$. At a maximizer $x_0\in J_0$,
$|F(x)|\ge A_0/2$ whenever
$|x-x_0|\le A_0/(2C_0B)$, by \eqref{eq:propM}.
For a small enough numerical constant $c>0$,
$J_1$ contains a one-sided or two-sided interval of length at least
$c\min\{1,A_0/B\}$ on which this lower bound holds. Because
$A_0\le C_0B$, integration over $O\supset J_1$ yields
\[
 \delta^2\ge cA_0^2\min\{1,A_0/B\}
 \ge c'\frac{A_0^3}{B}.
\]
Consequently
\begin{equation}\label{eq:smallonslit}
 \max_{J_0}|F|\le C_1B^{1/3}\delta^{2/3}.
\end{equation}

To propagate the smallness from $J_0$ to the reflected interval,
consider the slit domain $D=D_0\setminus J_0$. The bounded
subharmonic function $\log|F|$ (with value $-\infty$ at zeros)
obeys the usual two-constants inequality; equivalently one may apply
the maximum principle to $\log(|F|+\varepsilon)$ and let
$\varepsilon\downarrow0$ \cite{Vessella}. Let $\omega(z)$ denote the
harmonic measure of the two sides of the slit $J_0$ in $D$. On the
outer rectangle boundary, $|F|\le C_0B$; on either side of the slit,
$|F|\le A_0$. Thus, for $z\in D$,
\begin{equation}\label{eq:twoconstants}
 |F(z)|\le(C_0B)^{1-\omega(z)}A_0^{\omega(z)}.
\end{equation}
The harmonic measure $\omega$ is strictly positive at every interior
point of the connected bounded domain $D$, because the slit has
positive logarithmic capacity. By continuity and compactness, its
minimum on the part of $\overline{O^\sharp}$ outside the interior of
$J_0$ and inside $D$ is some $\omega_*>0$; on $J_0$ the stronger
bound \eqref{eq:smallonslit} already holds. Shrinking $\omega_*$ if
necessary and combining \eqref{eq:smallonslit} with
\eqref{eq:twoconstants} proves
\[
 \sup_{\overline{O^\sharp}}|F|
 \le C B^{1-\alpha}\delta^\alpha,
 \qquad \alpha=\frac23\omega_*\in(0,1).
\]
If $\delta=0$, the identity theorem gives $F\equiv0$ throughout $S$;
if $B=0$, then $u=0$. The remaining cases, including $\delta$ comparable
to $B$, follow by enlarging $C$ in \eqref{eq:propM}.

Finally \eqref{eq:reflection} and $\norm u_{L^1(I)}\le CB$ give
$0\le\mathcal E_d(u)\le CB^{2-\alpha}\delta^\alpha$.
In the positive atomic representation of \eqref{eq:reflection}, the
term at each $a_n$ contributes $b_ne^{-da_n}|U_n|^2\ge0$.
Keeping only these terms yields \eqref{eq:moments}. 
\end{proof}

\subsection{Quantitative recovery of the noninteger moments}
Put $a=4\pi/\beta$, $r=e^{-ax}$, $r_0=e^{-a\ell}$, and
\[
 q(r)=\begin{cases}u(-\log r/a)/(ar),&r_0<r<1,\\0,&0<r\le r_0.\end{cases}
\]
The zero trace of $u$ at both ends gives $q\in H^1(0,1)$ and
\begin{equation}\label{eq:qnorm}
 c\norm u_{L^2(I)}\le\norm q_2\le C\norm u_{L^2(I)},\qquad
 \norm q_{H^1(0,1)}\le C\norm u_H.
\end{equation}
Writing $\zeta_n=a_n/a$, the moments are
\[
 U_n=\int_0^1q(r)r^{\zeta_n}\dd r,\qquad
 \zeta_n=\sqrt{n^2+(m\beta/(2\pi))^2}+\frac{m\beta}{4\pi}.
\]
In particular $n+c_-\le\zeta_n\le n+c_+$ and $\zeta_{n+1}-\zeta_n\ge c_0>0$.

\begin{lemma}[A M\"untz approximation bound]\label{lem:muntzrate}
Let $E_N=\operatorname{span}\{r^{\zeta_1},\ldots,r^{\zeta_N}\}$ and let $\Pi_N$ be its orthogonal projection in $L^2(0,1)$. Then
\begin{equation}\label{eq:muntzrate}
 \norm{q-\Pi_Nq}_2\le\frac{C\norm q_{H^1(0,1)}}{\log(eN)}.
\end{equation}
The Gram matrix $G_N=(1/(\zeta_i+\zeta_j+1))_{i,j=1}^N$ satisfies
\begin{equation}\label{eq:gram}
 \norm{G_N^{-1}}_{\op}\le e^{CN^2\log(eN)},\qquad
 \norm{\Pi_Nq}_2^2=U^*G_N^{-1}U.
\end{equation}
\end{lemma}
\begin{proof}
The Cauchy determinant formula, applied to the Gram matrix after adjoining $r^k$, gives the exact distance identity
\begin{equation}\label{eq:muntzdist}
 \operatorname{dist}(r^k,E_N)
 =\frac1{\sqrt{2k+1}}\prod_{j=1}^N
 \frac{|k-\zeta_j|}{k+\zeta_j+1},\qquad k=0,1,\ldots.
\end{equation}
A coincident exponent gives zero and causes no difficulty. For $0\le k\le K$, retain factors with $j\ge J_K$, where $J_K\le C(K+1)$ is large enough that $\zeta_j>k$. The remaining factors obey
\[
 \frac{\zeta_j-k}{\zeta_j+k+1}
 \le\exp\left(-\frac{2k+1}{j+k+c_++1}\right).
\]
Comparison of the harmonic sum with its integral proves, when $N\ge C(K+1)$,
\[
 \operatorname{dist}(r^k,E_N)
 \le\left(\frac{C(K+1)}{N+1}\right)^{2k+1}
 \le\frac{C(K+1)}{N+1}.
\]

Let $p_K$ be the degree-$K$ shifted Legendre projection of $q$. The Legendre eigenvalue identity and integration by parts give
\[
 \norm{q-p_K}_2\le\frac{C\norm{q'}_2}{K+1}.
\]
Indeed the orthonormal coefficients $c_j$ satisfy
$\sum j(j+1)|c_j|^2\le\int_0^1r(1-r)|q'(r)|^2\dd r$.
The monomial coefficients of $p_K=\sum_{k=0}^KA_kr^k$ satisfy
$\sum|A_k|\le C_1^{K+1}\norm q_2$ for a fixed $C_1>1$, by the binomial formula for shifted Legendre polynomials. Hence
\[
 \operatorname{dist}(q,E_N)
 \le\frac{C\norm q_{H^1}}{K+1}
 +\frac{C(K+1)C_1^{K+1}}{N+1}\norm q_2.
\]
Choose $K=\lfloor\epsilon\log(N+1)\rfloor$ with $\epsilon$ small enough that $C_1^{K+1}\le C\sqrt N$. This proves \eqref{eq:muntzrate}. The distance identity \eqref{eq:muntzdist} and the preceding polynomial projection estimates give \eqref{eq:muntzrate}.

Finally
\[
 \det G_N=\frac{\prod_{i<j}(\zeta_j-\zeta_i)^2}
 {\prod_{i,j}(\zeta_i+\zeta_j+1)}.
\]
Uniform separation of the $\zeta_j$ and their linear growth give
$\det G_N\ge e^{-CN^2\log(eN)}$.
Using $\lambda_{\max}(G_N)\le\Tr G_N\le CN$ bounds its smallest eigenvalue from below by the same form after changing $C$. The projection identity in \eqref{eq:gram} follows by solving its normal equations.
\end{proof}

\begin{theorem}[Double-logarithmic Sobolev stability]\label{thm:doublelog}
Assume \eqref{eq:geom} and \eqref{eq:short}. For each $s>0$ there is $C_s$ such that, for all $u\in X_s$ and $M>0$ with $\norm u_{X_s}\le M$,
\begin{equation}\label{eq:doublelog}
 \norm u_H\le
 \frac{C_sM}{\LL(M/\norm{\T u}_Y)^{\theta_s}},\qquad
 \theta_s=\frac{s}{1+s},\qquad
 \LL(x)=1+\log(1+\log(1+x)).
\end{equation}
At zero data the right side is interpreted as zero. Constants depend on the fixed intervals, $m,\beta,L,s$, not on the source dimension or sign.
\end{theorem}
\begin{proof}
Initially assume only $\norm u_H\le B$. Equations~\eqref{eq:atoms} and \eqref{eq:moments} imply
\[
 \sum_{n\le N}|U_n|^2\le C(1+N)^2e^{CN}B^{2-\alpha}\delta^\alpha.
\]
Combining Lemma~\ref{lem:muntzrate} and \eqref{eq:qnorm} gives
\begin{equation}\label{eq:balance}
 \norm u_2\le\frac{CB}{\log(eN)}
 +\exp\{C_2N^2\log(eN)\}B^{1-\alpha/2}\delta^{\alpha/2}.
\end{equation}
For sufficiently small $\epsilon=\delta/B$, take
\[
 N=\left\lfloor c\sqrt{\frac{\log(1/\epsilon)}{\log\log(1/\epsilon)}}\right\rfloor
\]
with $c>0$ small. The second term is at most $CB\epsilon^{\alpha/4}$ and $\log(eN)\ge c'\log\log(1/\epsilon)$. Adjusting constants outside this range yields
$\norm u_2\le CB/\LL(B/\delta)$.

In the Dirichlet basis, H\"older's inequality gives the exact spectral interpolation
\[
 \norm u_H\le\norm u_2^{s/(1+s)}\norm u_{X_s}^{1/(1+s)}.
\]
Use $B=\mu_1^{-s/2}M$ and absorb fixed rescalings of the argument of $\LL$. If $\delta=0$, send $N\to\infty$ in \eqref{eq:balance} instead.
\end{proof}
Theorem~\ref{thm:doublelog} establishes the conditional bound under \eqref{eq:short}.

\begin{remark}[Origin of the double logarithm and possible improvement]\label{rem:rate-loss}
The propagation step \eqref{eq:propagate} converts the original data size $\delta$ into a \emph{power} $B^{1-\alpha}\delta^\alpha$ on an auxiliary interval. This step does not, by itself, force a double logarithm. The principal quantitative loss occurs in the subsequent reconstruction of the noninteger moments: our approximation estimate \eqref{eq:muntzrate} has only an inverse-$\log N$ rate, while \eqref{eq:gram} allows an amplification as large as $\exp(CN^2\log(eN))$. Balancing these two contributions in \eqref{eq:balance} selects $N$ of order $\sqrt{\log(B/\delta)/\log\log(B/\delta)}$ and yields the double logarithm. Sharper M\"untz--Jackson approximation or structured moment inversion might improve the rate, as suggested by the related mechanisms in \cite{GarciaFerreroRuland,Gendron,Choulli}; however, their estimates do not automatically apply to the present frequencies $\zeta_n$ together with the reflected thermal weights and Dirichlet-spectral prior. 
\end{remark}

\begin{corollary}[An explicit resolution function]\label{cor:explicitresolution}
Under the same geometry, the separation function in \eqref{eq:alpha} has the lower bound
\begin{equation}\label{eq:explicitalpha}
 \alpha_{M,s}(r)\ge
 \frac{2M}{\exp(\exp(q_r-1)-1)-1},\qquad
 q_r=\max\left\{2,\left(\frac{2C_sM}{r}\right)^{1/\theta_s}\right\}.
\end{equation}
Inserting this number in \eqref{eq:lambdar} and \eqref{eq:resolution} gives explicit sufficient coupling and information thresholds for every prescribed positive resolution.
\end{corollary}
\begin{proof}
Apply \eqref{eq:doublelog} to $u=v-w$, whose $X_s$ norm is at most $2M$, and solve
$\LL(2M/\norm{\T u})\le(2C_sM/r)^{1/\theta_s}$.
\end{proof}
The doubly exponential denominator in \eqref{eq:explicitalpha} leads to small sufficient thresholds; the constants follow from the stated analytic bounds.

\subsection{Continuous inversion and a spectral consequence}
For $t>0$ define
\[
 \gamma_s(t)=\inf_{u\in X_s,\norm u_H=1}
       (\norm{\T u}_Y^2+t\norm u_{X_s}^2).
\]
Theorem~\ref{thm:doublelog} implies, for small $t$,
\begin{equation}\label{eq:gammalower}
 \gamma_s(t)\ge c_st\LL(1/t)^{2\theta_s}.
\end{equation}
To verify this explicitly, a failure would give $t_j\downarrow0$, $\epsilon_j\downarrow0$ and unit $u_j$ with
$\delta_j^2+t_jq_j^2\le\epsilon_jt_j\LL(1/t_j)^{2\theta_s}$, where $q_j=\norm{u_j}_{X_s}\ge\mu_1^{s/2}$. Then $q_j\le\sqrt{\epsilon_j}\LL(1/t_j)^{\theta_s}$ and
$\LL(q_j/\delta_j)\ge c\LL(1/t_j)$. Applying \eqref{eq:doublelog} with prior size $q_j$ contradicts $\norm{u_j}_H=1$.

If $\norm{z^\delta-\T v_*}_Y\le\delta$ and $\norm{v_*}_{X_s}\le M$, the continuous Tikhonov minimizer
\[
 v_t^\delta=\arg\min_{v\in X_s}
 \left\{\tfrac12\norm{\T v-z^\delta}_Y^2+\tfrac t2\norm v_{X_s}^2\right\}
\]
exists uniquely. Its normal equation and Cauchy--Schwarz in the product energy norm give
\[
 \norm{v_t^\delta-v_*}_H^2\le\frac{\delta^2+tM^2}{\gamma_s(t)}.
\]
With $t=(\delta/M)^2$, this has the rate in \eqref{eq:doublelog}. The regularization argument is standard \cite{Engl}; the kernel-specific input is \eqref{eq:gammalower}.

\section{Fixed-coupling identifiability and its limitations}\label{sec:fixed}
\subsection{The exact residual problem}
For fixed $\lambda>0$ let $\widehat\nu_{\lambda,0}\propto e^{-\lambda(V_4+bV_2)}\mu$, and put
\[
 \Delta V_v=\int_{\mathscr T}\bigl(:(\phi+v)^4:_\beta-:\phi^4:_\beta\bigr).
\]
The calibrated shifted representation gives the exact local density
\begin{equation}\label{eq:conditional}
 \frac{\dd P_{\lambda,v}^O}{\dd P_{\lambda,0}^O}
 =\frac{\E_{\widehat\nu_{\lambda,0}}[e^{-\lambda\Delta V_v}\mid\OO]}
 {\E_{\widehat\nu_{\lambda,0}}e^{-\lambda\Delta V_v}}.
\end{equation}
This is a density on the equal-time \emph{commutative} field algebra, not a density matrix for the whole local quantum algebra. Equality of the two normalized conditional expectations is equivalent to equality of the corresponding measured laws. Showing that this implies $v=w$ at one fixed coupling would be sufficient for local quantum-state injectivity, but is not established here.

For smooth directions $u$, differentiating \eqref{eq:conditional} gives the measured Fisher information
\begin{equation}\label{eq:fisher}
 \mathcal I_{\lambda,v}^O(u,u)
 =\lambda^2\Var_{\widehat\nu_{\lambda,v}}
       \bigl(\E_{\widehat\nu_{\lambda,v}}[Q_v[u]\mid\OO]\bigr),
 \qquad Q_v[u]=4\int_{\mathscr T}u:(\phi+v)^3:_\beta.
\end{equation}
The conditional-score identity is standard \cite{Louis}. Its integrability here follows from Lemma~\ref{lem:Nelson}. Positivity of the conditional variance in every nonzero source direction is an open problem. The term $\lambda^2\Var(Q_v[u])$ in the parameter-dependent shifted full model generally differs from the physical BKM metric; the local observation law is invariant under that shift.

\subsection{A limitation of perturbative inference}
\begin{proposition}[An analytic remainder can destroy injectivity]\label{prop:counter}
There are bounded compact operators $T,R$ on $\ell^2$, with $T$ injective, such that
$M_\lambda=-24\lambda T+24\lambda^2R$ has a nonzero kernel at arbitrarily small positive couplings, despite a uniform $O(\lambda^2)$ derivative remainder bound.
\end{proposition}
\begin{proof}
Take $Te_n=e^{-n}e_n$ and $Re_n=ne^{-n}e_n$. Both are bounded and compact, and $M_{1/n}e_n=0$.
\end{proof}
The example isolates the additional structure required to deduce fixed-coupling infinite-dimensional uniqueness from Theorem~\ref{thm:response} and injectivity of $\T$. Similarly, the diagonal analytic kernels with singular values $e^{-n}$ and $e^{-n^2}$ both have exponential upper bounds but different logarithmic inverse scales. An upper bound does not select a matching stability rate.

\subsection{A sufficient structural condition}
For a fixed coupling and reference source $v_0$, let
$\eta(v)=\sqrt{\dd P_{\lambda,v}^O/\dd P_{\lambda,v_0}^O}$ and $\mathcal A_0=D\eta(v_0)$, when these derivatives exist. If on a convex source set
\begin{equation}\label{eq:alignment}
 \ip{D\eta(v)[u]}{\mathcal A_0u}\ge c\norm{\mathcal A_0u}_2^2\quad(c>0),
\end{equation}
then line integration gives
\[
 H^2(P_{\lambda,v}^O,P_{\lambda,w}^O)
 \ge c^2\norm{\mathcal A_0(v-w)}_2^2
 =\tfrac{c^2}{4}\mathcal I_{\lambda,v_0}^O(v-w,v-w).
\]
Indeed pair $\eta(v)-\eta(w)$ with $\mathcal A_0(v-w)$ and apply Cauchy--Schwarz. Whether \eqref{eq:alignment} holds for the interacting field remains open.

\begin{problem}\label{prob:fixed}
For a fixed positive coupling, is the conditional law in \eqref{eq:conditional} injective on an exterior Sobolev prior? If so, what uniform conditional stability modulus can be proved? Is there a useful field-specific estimate such as \eqref{eq:alignment}, or a different completeness mechanism for the conditional expectation?
\end{problem}

\section{Conclusion}
Quadratic-reference calibration isolates a cubic local response from the thermal mass shift. Its continuum weak-coupling expansion, kernel injectivity and analytic smoothing yield local relative-entropy separation for finite-dimensional source spaces and for prescribed Sobolev resolutions. A reflected positive energy gives conditional stability under the additional interval-length assumption. Fixed-coupling infinite-dimensional identifiability is formulated in Problem~\ref{prob:fixed}.

\section*{Acknowledgments and funding}
Rubin Zhang acknowledges support from the China Scholarship Council (CSC). Ning Zhang acknowledges support from the National Natural Science Foundation of China (Grant No.~62535007) and Huazhong University of Science and Technology (Grant No.~5003011074).

\section*{Declaration of AI assistance}
During manuscript preparation, ChatGPT (OpenAI) assisted with literature review, drafting and editing, organization and checking of mathematical derivations, and software and figure generation. The authors take responsibility for the content and accuracy of the manuscript and accompanying materials.

\section*{Data and code availability}
No external datasets were used. The Python script \texttt{reproduce\_figure.py}, the CSV values \texttt{figure1\_data.csv}, and reproduction instructions are supplied as supplementary material (and as ancillary files in the \texttt{anc/} directory of the arXiv source). The vector Figure~\ref{fig:singvalues} is supplied as \texttt{bkm\_singular\_values.pdf}.

\end{document}